\documentclass[a4paper]{amsart}

\usepackage[latin1]{inputenc}
\usepackage[T1]{fontenc}
\usepackage[linktocpage]{hyperref} 
	\hypersetup{colorlinks=true,linkcolor=blue,citecolor=blue,filecolor=magenta,urlcolor=blue}    
\usepackage{geometry}       
\usepackage[english]{babel}  
\usepackage{amsmath,amsfonts,amssymb,amsthm}
\usepackage{mathrsfs}
\usepackage{dsfont}
\usepackage{tikz,tkz-graph,tikz-cd}

\newtheorem{thm}{Theorem}[section]

\newtheorem{propo}[thm]{Proposition}

\newtheorem{cor}[thm]{Corollary}

\theoremstyle{definition}
\newtheorem{de}[thm]{Definition}

\theoremstyle{remark}
\newtheorem{rmk}[thm]{Remark}

\newcommand{\CC}{\mathds{C}}

\newcommand{\ZZ}{\mathds{Z}}

\newcommand{\FF}{\mathds{F}}
\newcommand{\PP}{\mathds{P}}

\newcommand{\BB}{\mathds{B}}
\newcommand{\A}{\mathcal{A}}
\renewcommand{\P}{\mathcal{P}}
\newcommand{\Q}{\mathcal{Q}}

\newcommand{\D}{\mathcal{D}}
\newcommand{\C}{\mathcal{C}}

\newcommand{\I}{\mathcal{I}}
\newcommand{\J}{\mathcal{J}}

\renewcommand{\S}{\mathcal{S}}

\renewcommand{\L}{\mathcal{L}}
\newcommand{\inv}{^{-1}}

\renewcommand{\epsilon}{\varepsilon}

\DeclareMathOperator{\HH}{H}

\DeclareMathOperator{\lk}{lk}

\DeclareMathOperator{\Tub}{Tub}
\DeclareMathOperator{\Irr}{Irr}
\DeclareMathOperator{\Top}{top}
\DeclareMathOperator{\Sing}{Sing}

\begin{document}

\title{A linking invariant for algebraic curves}

\author[B. Guerville-Ball\'e]{Beno\^it Guerville-Ball\'e}
         \address{Instituto de Ci\^encias Matem\'aticas e de Computa\c{c}\~ao
						Universidade de S\~ao Paulo
						Avenida Trabalhador Sancarlense, 400 - Centro
						S\~ao Carlos - SP, 13566-590 (Brazil)}
         \email{benoit.guerville-balle@math.cnrs.fr }
\author[J.B. Meilhan]{Jean-Baptiste Meilhan}
         \address{Universit\'e Grenoble Alpes, IF, 38000 Grenoble, France}
         \email{jean-baptiste.meilhan@univ-grenoble-alpes.fr}
        
        \thanks{The first author was supported by a JSPS Postdoctoral Fellowship and is currently supported by the postdoctoral grant \#2017/15369-0 of the Funda\c{c}\~ao de Amparo \`a Pesquisa do Estado de S\~ao Paulo (FAPESP). The second author was supported by the French ANR research project ``VasKho'' ANR-11-JS01-00201.}				


\selectlanguage{english}
\begin{abstract}
We construct a topological invariant of algebraic plane curves, which is in some sense an adaptation of the linking number of knot theory. This invariant is shown to be a generalization of the $\I$-invariant of line arrangements developed by the first author with Artal and Florens. We give two practical tools for computing this invariant, using a modification of the usual braid monodromy or using the connected numbers introduced by Shirane. As an application, we show that this invariant distinguishes several Zariski pairs, i.e. pairs of curves having same combinatorics, yet different topologies. The former is the well known Zariski pair found by Artal, composed of a smooth cubic with 3 tangent lines at its inflexion points. The latter is formed by a smooth quartic and 3 bitangents.
\end{abstract}

\maketitle

\section*{Introduction}
	The topological study of algebraic plane curves was initiated at the beginning of the 20th century by Klein and Poincar\'e. One of the main questions is to understand the relationship between the combinatorics and the topology of a curve. It is known, since the seminal work of Zariski~\cite{zariski,zariski1,zariski2}, that the topological type of the embedding of an algebraic curve in the complex projective plane is not determined by the combinatorics. Indeed, Zariski constructed two sextics with $6$ cusps having same combinatorics, and proved that the fundamental group of their complements are not isomorphic. Geometrically, these two curves are distinguished by the fact that the cusps in the first curve lie on a conic, while they do not in the second curve. Since this historical example, using various methods, numerous examples of pairs of algebraic curves having same combinatorics but different topologies have been found, see for example Artal, Cogolludo and Tokunaga~\cite{survey}, Cassou-Nogu\`es, Eyral and Oka~\cite{CasEyrOka_Topology}, Degtyarev~\cite{Deg_Deformations}, Oka~\cite{Oka_Transforms}, Shimada~\cite{Shimada_Fundamental}, or the first author \cite{guervilleGT}. Artal suggests in~\cite{artal} to call such examples \emph{Zariski pairs}.\\

The topology of curves in $\CC\PP^2$ is intimately connected to the topology of knots and links in $S^3$. Several tools are indeed shared by these two domains, such as the homology or the fundamental group of the complement, the Alexander polynomial or module, although they usually have rather different behaviours.

Recently, Artal, Florens and the first author defined a topological invariant of line arrangements (i.e. algebraic plane curves with only irreducible components of degree $1$) which is in some sense modelled on the linking number of knot theory~\cite{invariant}. This invariant was then successfully used in~\cite{guervilleGT} to distinguish a new Zariski pair of line arrangements. In the present paper, we construct another invariant adapting the linking number to the more general case of algebraic plane curves. In the case of a line arrangement, this invariant is shown to be equivalent to the invariant of~\cite{invariant}, thus providing a generalization of this earlier work through a different adaptation of the linking number.

The construction of our linking invariant can be roughly outlined as follows. Consider a reducible algebraic curve decomposed in two subcurves $\C$ and $\D$, and pick a topological cycle $\gamma$ in the subcurve $\C$. The basic idea is to consider the image of a certain coset of $\gamma$ in the first homology group of $\CC\PP^2\setminus\D$. More precisely, this set is regarded in the quotient of $\HH_1(\CC\PP^2\setminus\D)$ by an appropriate \emph{indeterminacy subgroup}~$\J_\C$, which controls the topological differences among the various cycles in the considered coset of~$\gamma$. This define the \emph{linking invariant of $\C$ with $\D$ along $\gamma$}, which is an invariant of the pair $(\CC\PP^2,\C\cup \D)$.

Our construction thus builds on a rather elementary idea, and is not technically involved. Remarkable is rather the fact that it reveals quite efficient in practice, despite its apparent simplicity. We mention below several applications of the linking invariant on concrete examples of Zariski pairs of various natures. 

This linking invariant has a nice behaviour for some particular choices of curve or cycle. In the case of rational curves, that is, curves whose irreducible components have all genus $0$, the linking invariant is indeed a single homology class rather than a coset. This allows us to prove the equivalence with the $\I$-invariant of~\cite{invariant} in the case of line arrangements.\\

From a practical viewpoint, we provide two methods of computation of this linking invariant. The first one is based on a topological construction using an adaptation of the braid monodromy. This makes a concrete connection between our invariant and the usual linking number of knot theory. The second method is algebraic and comes from the relation, observed in~\cite{GueShi}, between our linking invariant, the connected numbers and the splitting numbers introduced by Shirane in~\cite{shirane,shirane_connected}.

To illustrate the efficiency of this adaptation of the linking number to algebraic curves, we use it to distinguish two examples of Zariski pairs. The first example is formed by the well known $3$-Artal curves introduced by Artal in~\cite{artal}. They are composed of a smooth cubic and three inflexional tangent; in the first curve the considered inflexion points are collinear, while they are not in the second one. The computation of the linking of the three lines with the cubic is made using the above mentioned algebraic method. The second example of Zariski pair is formed by a smooth quartic and three bitangents. 
These curves have been very recently studied in~\cite{BTY}. For that example, we use the topological method based on the linking number of knot theory.\\

After an earlier version of this paper was circulated, our linking invariant (then called \emph{linking set}) has been further studied, and being used to distinguish other examples of Zariski pairs. 

In~\cite{shirane} Shirane introduces the \emph{splitting numbers} and detects the $\pi_1$-equivalent Zariski $k$-plets suggested by Shimada~\cite{Shimada_Equisingular}. By proving that the splitting numbers and the linking invariant are equivalent (in some particular cases), the first author and Shirane obtain in~\cite{GueShi} that the linking invariant distinguishes the Shirane-Shimada $\pi_1$-equivalent Zariski $k$-plets. This implies that the linking invariant is not determined by the fundamental group of the complement. 

The linking invariant is also used in~\cite{BGST} to classify the topology of the $k$-Artal curves (i.e. a smooth cubic and $k$ inflectional tangent lines). Furthermore, Shirane constructed recently in~\cite{shirane_connected} an adaptation of the splitting number, called the \emph{connected numbers}, which allows to classify the topology of the Artal curves of degree $b$ (i.e. smooth curves of degree $b$ and with three total inflectional tangent lines). Here again, the proofs of~\cite{GueShi} imply that the linking invariant can distinguish the Artal curves of degree $b$. 

In the particular case of line arrangements, the linking invariant (in the form of the $\I$-invariant) has been successfully used in~\cite{guervilleGT} to detect a Zariski pair of $12$ lines. Recently, the first author and Viu-Sos gave an effective diagrammatic reformulation of this invariant in the particular case of real line arrangements, see~\cite{GueViu}. Using this reformulation, they provide $10$ examples of complexified real Zariski pairs.\\

\noindent\textbf{Convention.}\mbox{}
All homology groups are to be understood with integral coefficients, and this will be omitted in the notation.

\section{The linking invariant}\label{sec:linking}
	\subsection{Preliminaries}\label{sec:preli}

Let $\C$ be an algebraic plane curves, possibly non-reduced. Following~\cite{survey}, we define the combinatorics of $\C$ as the data:
\begin{equation*}
	\left( \Irr(\C), \deg, \Sing(\C), \Sigma_{\Top}, \sigma_{\Top}, \{\C(P)\}_{P\in \Sing(\C)}, \{\beta_P\}_{P\in \Sing(\C)} \right), 
\end{equation*}
where:
\begin{itemize}
	\item $\Irr(\C)$ is the set of all irreducible components of $\C$,  
	\item $\deg$ assigns to each irreducible component its degree, 
	\item $\Sing(\C)$ is the set of all singular points of $\C$, 
	\item $\Sigma_{\Top}$ is the set of topological types of singular points of $\C$, and $\sigma_{\Top}$ assigns to each singular point its topological type, 
	\item for each singular point $P$ of $\C$, $\C(P)$ is the set of local branches of $\C$ at $P$, and	$\beta_P$ assigns to each local branch at $P$ the global irreducible component containing it.
\end{itemize}
Two curves have the same combinatorics if there exist bijections between their sets $\Irr$ and $\Sing$ of irreducible components and singular points, which are compatible with the sets $\C(P)$ and $\Sigma_{\Top}$ of local branches and topological types, and with the assignments $\deg$, $\sigma_{\Top}$ and $\{\beta_P\}$ in the natural way; see \cite[Rem.~3]{survey} for details.

We also associate to the curve $\C$ the \emph{intersection graph} $\Gamma_\C$ of its irreducible components. This is a bipartite graph whose first set of vertices, called \emph{component-vertices}, corresponds to the irreducible components of $\C$, while the second set of vertices, called \emph{point-vertices}, corresponds to the singular points of $\C$ contained in two distinct irreducible components. An edge of $\Gamma_\C$ joins a point-vertex and a component-vertex if and only if the corresponding singular point is contained in the corresponding irreducible component. Note that the information encoded in $\Gamma_\C$ are contained (but not equivalent) to the combinatorics of $\C$. For example, the information given by $\sigma_{\Top}$ is not contained in $\Gamma_\C$.

A \emph{cycle} of $\Gamma_\C$ is a (non necessarily connected) closed oriented walk without repeated edges. Note that this includes the case of a single vertex. A (combinatorial) cycle of $\Gamma_\C$ can be lifted to a (topological) cycle on the curve $\C$, \emph{i.e.} an oriented closed loop in $\C$, although it is not uniquely determined in general. In particular, a topological lift of a combinatorial cycle has a natural induced orientation only if the combinatorial cycle is not simply connected. 

In what follows, we will be mainly interested in reducible algebraic curves $\C\cup \D$, which decompose into two subcurves $\C$ and $\D$ (without common irreducible component). In this context, a cycle of $\Gamma_\C$ is simply a cycle lying in the subgraph $\Gamma_\C$ of $\Gamma_{\C\cup \D}$. On one hand, such a cycle is called \emph{maximal} if it contains all component-vertices of $\Gamma_\C$; in other words, a maximal cycle in $\Gamma_\C$ lifts to a topological cycle in $\C\cup \D$ that intersects the smooth part of all irreducible components of $\C$. On the other hand, a cycle of $\Gamma_{\C\cup \D}$ is said to \emph{avoid $\D$} if it lies in $\Gamma_\C$ (and is thus disjoint from all component-vertices of $\D$) \emph{and} avoids all point-vertices of $\C\cap \D$. 

In this paper, by a homeomorphism $\phi$ between two such reducible algebraic curves $\C_1\cup \D_1$ and $\C_2\cup \D_2$, we will always mean an ambient homeomorphism of $\CC\PP^2$ which sends $\C_1$ to $\C_2$ and $\D_1$ to $\D_2$. Furthermore, we will denote by $\phi_{\Gamma}: \Gamma_{\C_1\cup \D_1}\rightarrow \Gamma_{\C_2\cup \D_2}$ the induced map at the combinatorial level. Note that, if $\phi$ is orientation preserving, then $\phi_{\Gamma}$ preserves the cycle orientation. 

\subsection{The linking invariant}

Let $\C\cup \D$ be a reducible algebraic curve, decomposed into two subcurves $\C$ and $\D$.

Consider the inclusion maps $i:\C\setminus\D\hookrightarrow\C$ and $j:\C\setminus\D\hookrightarrow\CC\PP^2\setminus\D$, and denote respectively by $i_*$ and $j_*$ the induced map on the first homology groups. 
Note that $\ker(i_*)$ identifies with $\HH_1(\partial \bigcup_{C\in\Irr(\C)} C\setminus \D)\simeq \bigoplus_{C\in\Irr(\C)} \HH_1\left(\partial (C\setminus \D)\right)$ in $\HH_1(\C\setminus\D)$. 
\begin{de}
 The \emph{indeterminacy subgroup} with respect to $\C$, denoted by $\J_\C$, is the subgroup of $\HH_1(\CC\PP^2\setminus\D)$ defined as the image of $\bigoplus_{C\in\Irr(\C)} \HH_1\left(\partial (C\setminus \D)\right)$  by $j_*$. 
\end{de}

Now, let $\gamma$ be a maximal cycle in $\Gamma_\C$ avoiding $\D$. Pick a topological lift $\widetilde{\gamma}$ of $\gamma$ on the curve $\C\subset \C\cup \D$. By assumption, $\widetilde{\gamma}$ lies in $\C\setminus \D$, and intersects the smooth part of all irreducible components of $\C$. 

For brevity, we simply denote by $[\gamma]$ the image of $\widetilde{\gamma}$ in $\HH_1(\CC\PP^2\setminus\D) / \J_\C$. We also denote by $\I_\C$ the image of $\bigoplus_{C\in\Irr(\C)} \HH_1(C\setminus \D)$ by $j_*$, composed with the projection map $\HH_1(\CC\PP^2\setminus\D)\rightarrow \HH_1(\CC\PP^2\setminus\D)/\J_\C$. 

\begin{de}\label{def:linking}
	The \emph{oriented linking of $\C$ with $\D$ along $\gamma$}, denoted by $\lk_\gamma(\C,\D)$, is the coset of $\I_\C$ in $\HH_1(\CC\PP^2\setminus\D) / \J_\C$ with respect to $[\gamma]$. In other words, 
	\begin{equation*}
		\lk_\gamma(\C,\D) = [\gamma]\I_\C\ \in\ \HH_1(\CC\PP^2\setminus\D) / \J_\C.
	\end{equation*}
\end{de}

\begin{thm}\label{thm:well}
The above formula is well-defined, \emph{i.e.} $\lk_\gamma(\C,\D)$ does not depend on the choice of topological lift of $\gamma$.
\end{thm}
\begin{proof}
Let $\widetilde{\gamma}_1$ and $\widetilde{\gamma}_2$ be two topological lifts of $\gamma$, and let $[\gamma]_1$ and $[\gamma]_2$ denote their homology classes in  $\HH_1(\CC\PP^2\setminus\D)$. There are essentially two ways in which $[\gamma]_1$ and $[\gamma]_2$ may differ. If $\widetilde{\gamma}_1$ and $\widetilde{\gamma}_2$ have same homology class in $\HH_1(\C)$, then they differ by elements of $\cup_{C\in\Irr(\C)} \partial (C\setminus \D)$, so that $[\gamma]_1$ and $[\gamma]_2$ differ by an element of the indeterminacy subgroup $\J_C$. Now, if $\widetilde{\gamma}_1$ and $\widetilde{\gamma}_2$ have different homology classes in $\HH_1(\C)$, then the difference is mapped in $\I_\C$ by $j_\ast$, so that  $\widetilde{\gamma}_1$ and $\widetilde{\gamma}_2$ yield the same coset of $\I_\C$ in $\HH_1(\CC\PP^2\setminus\D) / \J_C$.
\end{proof}

\begin{rmk}\label{rem:stress}
We stress that the neither of the two assumptions made here, that $\gamma$ is maximal and that it avoids $\D$, is necessary to define our invariant -- this is discussed in Remark \ref{rem:maxi} and in Section  \ref{sec:general} below. But, on one hand, these assumptions turn out to greatly simplify the exposition and, on the other hand, all the relevant topological information on $\C\cup \D$ are already essentially detected by this simple version of our invariant. As a matter of fact, all the examples of this paper will involve the above assumptions. 
\end{rmk}

We have the following description of $\J_\C$. 
\begin{propo}\label{propo:indeterminacy}
	The indeterminacy subgroup $\J_\C$ is spanned by the elements of the form:
	\begin{equation*}
		\sum_{d\in \D(P)} I_P(b,d) \cdot m_{\beta_P(d)}, \quad \text{for all $P\in\C\cap\D$ and all $b\in \C(P)$,}
	\end{equation*}
	where 
	$I_P(b,d)$ denotes the intersection multiplicity of the local branches $b$ and $d$ at $P$, and $m_{\beta_P(d)}$ is given by  a meridian of the irreducible component $\beta_P(d)$ of $\D$ containing $d$. 
\end{propo}

\begin{proof}
The indeterminacy subgroup is the image of  $\bigoplus_{C\in\Irr(\C)} \HH_1\left(\partial (C\setminus \D)\right)$ in $\HH_1(\CC\PP^2\setminus\D)$. It is thus generated by the class of the cycles in $C\in\Irr(\C)$ around the points $P\in\C\cap\D$. Pick such a singular point $P$, and consider a small sphere $S$ around $P$. Each local branch $b$ of $\C$ at $P$ intersects $S$ along a knot $K_b$, and it is well-known that, for each local branch $d$ of $\D$ at $P$, the intersection of $b\cup d$ with $S$ is an oriented two-component link whose linking number is precisely $I(b,d)$ (see \cite[pp.~439]{brieskorn}). Hence the homology class of the knot $K_b$ in $\HH_1(\CC\PP^2\setminus\D)$ is given by $\sum_{d\in \D(P)} I(b,d) \cdot m_{\beta_P(d)}$, and the result follows. 
\end{proof}
 
As a consequence of Proposition \ref{propo:indeterminacy}, the group $\HH_1(\CC\PP^2\setminus\D) / \J_\C$ is determined by the combinatorics of the curve $\C\cup\D$. So we can use the linking invariant to compare the topology of curves with the same combinatorics. Indeed, we have the following theorem, which implies that the linking invariant is an invariant of the oriented topology of $(\CC\PP^2,\C\cup \D)$.  
\begin{thm}\label{thm:invariance}
Let $\phi$ be an orientation-preserving homeomorphism between two algebraic curves $\C_1\cup \D_1$ and $\C_2\cup \D_2$. Then  $\phi$ induces an isomorphism $\phi_*$ between $\HH_1(\CC\PP^2\setminus\D_1)$ and $\HH_1(\CC\PP^2\setminus\D_2)$ mapping $\J_{\C_1}$ to $\J_{\C_2}$, and for any cycle $\gamma_1\in\Gamma_{\C_1}$ avoiding $\D_1$, we have 
	\begin{equation*}
		\overline{\phi_*}(\lk_{\gamma_1}(\C_1,\D_1))=\lk_{\phi_{\Gamma}(\gamma_1)}(\C_2,\D_2),
	\end{equation*}
	where $\overline{\phi_*}$ is map induced by $\phi_*$ on the quotients by the indeterminacy subgroups.
\end{thm}

\begin{proof}
By definition, the homeomorphism $\phi:\CC\PP^2\rightarrow\CC\PP^2$ maps $\D_1$ to $\D_2$, so it induces an isomorphism $\phi_*$ between $\HH_1(\CC\PP^2\setminus\D_1)$ and $\HH_1(\CC\PP^2\setminus\D_2)$. Furthermore, for each $C_1\in \Irr(\C_1)$ with image $C_2=\phi(C_1)\in \Irr(\C_2)$, we have that $\phi$ maps $C_1\cap \D_1$ to $C_2\cap \D_2$, and maps $\partial (C_1\setminus \D_1)$ to $\partial (C_2\setminus \D_2)$; this implies that $\phi_*$ maps $\J_{\C_1}$ to $\J_{\C_2}$. Now, $\phi$ maps any (oriented) lift of $\gamma_1$ to a cycle on $\C_2$ which is a lift of $\phi_{\Gamma}(\gamma_1)$, respecting the orientation. Since the linking invariant does not depend on the choice of lift, the result follows. 
\end{proof}

\begin{rmk}\label{rem:maxi}
As mentioned in Remark \ref{rem:stress}, the cycle $\gamma$ in Definition \ref{def:linking} does not need be maximal. Indeed, if $\Gamma$ does not contain all component-vertices of $\Gamma_\C$, then 
the coset $[\gamma]\I_\C$ still yields an invariant of the oriented topology of $(\CC\PP^2,\C\cup \D)$. But in this case a finer invariant is given by regarding the curve $\C\cup \D$ as decomposed into the union of $\C_\gamma$ and $(\C\cup D)\setminus \C_\gamma$, where $\C_\gamma$ denotes the union of all irreducible components of $\C$ intersecting $\gamma$. 
\end{rmk} 

\begin{rmk}\label{rem:non-oriented}
The linking invariant of Definition \ref{def:linking} is an invariant of the oriented topology of $(\CC\PP^2,\C\cup \D)$. If $\gamma$ is simply connected, however, the linking $\lk_\gamma(\C,\D)$ of $\C$ with $\D$ along $\gamma$ is a topological invariant of $(\CC\PP^2,\C\cup\D)$, since any choice of orientation of a topological lift yields the same coset. In general, we can easily remove the condition of orientation, simply by considering 
	\begin{equation*}
		-[\gamma]\I_\C \cup [\gamma]\I_\C\ \in\ \HH_1(\CC\PP^2\setminus\D) / \J_\C, 
	\end{equation*}
which is clearly an invariant that doesn't depend on the orientation of $\gamma$, but only on its combinatorics. As a corollary to Theorem~\ref{thm:invariance}, this \emph{non-oriented linking} is a topological invariant of the pair $(\CC\PP^2,\C\cup \D)$.
\end{rmk}

\subsection{Two variants}
We now discuss two variants of our linking invariant.
The first one is a `global' version which doesn't rely upon the choice of a cycle; 
the second one is a generalization, where we allow arbitrary cycles.

As in the previous section, $\C\cup \D$ will denote here a reducible algebraic curve decomposed into two subcurves $\C$ and $\D$.
\subsubsection{Global linking}
We can define the following coarser invariant, which is a `global' version of the linking invariant, in the sense that in doesn't involve the choice of a cycle.
Recall that $j_\ast$ is the map induced in homology by the inclusion of 
$\C\setminus \D$ in $\CC\PP^2\setminus \D$. 
\begin{de}
 The \emph{global linking} of $\C$ with $\D$, denoted by $\L(\C,\D)$, 
 is the class of Im$j_\ast$ in $\HH_1(\CC\PP^2\setminus\D)/\J_\C$. 
\end{de}
Remark that the global linking of $\C$ with $\D$ can also be defined as the union, over all cycles $\gamma$ in $\Gamma_\C$ avoiding $\D$, of the linking invariants of $\C$ with $\D$ along $\gamma$ . 

The invariance of the global linking is a direct consequence of the proof of Theorem \ref{thm:invariance}:
\begin{thm}
	Let $\phi$ be a homeomorphism between two curves $\C_1\cup \D_1$ and $\C_2\cup \D_2$. We have
	\begin{equation*}
		\overline{\phi_*}(\L(\C_1,\D_1))=\L(\C_2,\D_2),
	\end{equation*}
	where $\overline{\phi_*}$ is the map induced by $\phi_*$ on the quotients by the indeterminacy subgroups.

	Hence the global linking of $\C$ with $\D$ is a topological invariant of the pair $(\CC\PP^2,\C\cup \D)$.
\end{thm}
\begin{rmk}\label{rem:particular cases}
Notice that if $\C$ is an irreducible curve, then $\lk_*(\C,\D)=\L(\C,\D)$, where $*$ denotes the unique component-vertex of $\Gamma_\C$.
\end{rmk}

\subsubsection{Linking along an arbitrary cycle}\label{sec:general}
In the above definition of the linking invariant, we assumed throughout that the cycle $\gamma$ avoids all point-vertices of $\Gamma_\C$ corresponding to singularities in $\C\cap \D$. 
Although this will not be used in the main examples of this paper, we outline here how the construction can be easily generalized to arbitrary cycles.

Let $\gamma$ be any cycle in $\Gamma_\C$.
Denote by $S_\gamma$ the set of singularities in $\C\cap \D$ whose corresponding point-vertices are contained in $\gamma$. 
We define the subgroup $\J_\gamma$ of $\HH_1(\CC\PP^2\setminus\D)$ as 
$$ \J_\gamma = \langle m_{\beta_P(d)}\text{, for all $P\in\S_{\gamma}$ and all $d\in \D(P)$} \rangle. $$
\noindent (here we make use of the same notation as in Proposition \ref{propo:indeterminacy}.)
\begin{de}
	The \emph{$\gamma$-indeterminacy subgroup} is the subgroup $\J_{(\C,\gamma)}$ of $\HH_1(\CC\PP^2\setminus\D)$ generated by $\J_\C\cup \J_\gamma$. 
\end{de}

Now, we need a slightly generalized notion of topological lift for the cycle $\gamma$. 
Specifically, for each singular point $P$ in $\S_{\gamma}$, pick a small closed $4$-ball $\BB_P$ centered at $P$. 
A \emph{$\D$-avoiding lift} of $\gamma$ is a cycle in $\CC\PP^2\setminus \D$ which coincides with a topological lift outside $\cup_{P\in \S_\gamma} \BB_P$ (and in particular lies in $\C$), 
and whose intersection with each $4$-ball $\BB_P$ is an arc lying in the boundary of $\BB_P$. 
So, roughly speaking, such a cycle differs from a topological lift of $\gamma$ by locally pushing it away from the curve $\C\cup \D$, so that it avoids the singularities in $\S_\gamma$. 

Now, using this refined indeterminacy subgroup and generalized notion of lift, the exact same construction yields an invariant: 
the (oriented) linking of $\C$ with $\D$ along an arbitrary cycle $\gamma$ is the coset 
	\begin{equation*}
		\lk_\gamma(\C,\D) = [\gamma]\I_\C\ \in\ \HH_1(\CC\PP^2\setminus\D) / \J_{(\C,\gamma)}, 
	\end{equation*}
where $[\gamma]$ denotes the image of a $\D$-avoiding lift of $\gamma$ in $\HH_1(\CC\PP^2\setminus \D) / \J_{(\C,\gamma)}$, 
and where $\I_\C$ is the image of $\bigoplus_{C\in\Irr(\C)} \HH_1(C\setminus \D)$ by $j_*$, seen in $\HH_1(\CC\PP^2\setminus\D)/\J_{(\C,\gamma)}$. 

The proof that this is well-defined, \emph{i.e.} does not depend on the choice of the $\D$-avoiding lift of $\gamma$, is completely similar to the proof of Theorem \ref{thm:well}, and readily follows from the definition of $\J_{(\C,\gamma)}$. 
The only difference here is that, if $\gamma$ passes through a point-vertex $P$ in $\C\cap \D$, then two $\D$-avoiding lifts of $\gamma$ may only differ by a copy of a meridian $m_{\beta_P(d)}$, for any local branch $d$ in $\D(P)$. 
Considering these cycles in $\HH_1(\CC\PP^2\setminus\D)$, we have by definition that the difference lies precisely in $\J_\gamma$, hence in the $\gamma$-indeterminacy subgroup $\J_{(\C,\gamma)}$. 
\begin{rmk}\label{rem:yob}
In the case where $\gamma$ avoids all point-vertices in $\C\cap \D$, then the $\gamma$-indeterminacy subgroup $\J_{(\C,\gamma)}$ coincides with the original indeterminacy subgroup $\J_\C$, and we recover the invariant of Definition \ref{def:linking}. 
\end{rmk}

\subsection{Line arrangements and the $\I$-invariant}

When the curve $\C$ is a rational curve, that is when all the irreducible components of $\C$ have genus zero, the set $\I_\C$ is trivial. The linking of $\C$ and $\D$ along a cycle $\gamma$ is thus the class $[\gamma]$ of a lift  of $\gamma$ in $\HH_1(\CC\PP^2\setminus\D) / \J_\C$. 

This is, in particular, the case when $\C\cup\D$ is a line arrangement, \emph{i.e.} when all the irreducible components of $\C$ and $\D$ are of degree 1. In this particular case, another linking invariant, called the $\I$-invariant, has been defined by Artal, Florens and the first author in~\cite{invariant}. We will prove here that the present linking invariant generalizes the $\I$-invariant. 

First, let us recall some terminologies introduced in~\cite{invariant} to define the $\I$-invariant.
Let $\A$ be a line arrangement. 
Recall that $\HH_1(\CC\PP^2\setminus\A)$ is free of rank $\vert \Irr(\A)\vert-1$ and is generated by the set of all meridiens $m_L$ for $L\in\Irr(\A)$.
\begin{de}\label{def:inner-cyclic}
	Let $\A$ be a line arrangement, $\xi:\HH_1(\CC\PP^2\setminus\A)\rightarrow \CC^*$ be a non-trivial character and $\gamma$ be a cycle of the intersection graph $\Gamma_\A$. 
	The triple $(\A,\xi,\gamma)$ is an inner-cyclic arrangement if	
	\begin{enumerate}
		\item For each singular point $P$ of $\A$ with associated point-vertex in $\gamma$, 
		
		$\quad \xi(m_\ell)=1$, for any line $\ell$ of $\A$ containg $P$, 
		\item For each line $L$ of $\A$ with associated component-vertex in $\gamma$, 
		
		\begin{enumerate}
		\item[(i)] $\xi(m_L)=1$, 		
		\item[(ii)] $\prod\limits_{\ell \ni S} \xi(m_\ell)=1$, for any singular point $S$ in $L$, 
		where the product runs over all the lines $\ell$ in $\A$ containing $S$.	
		\end{enumerate}
	\end{enumerate}
\end{de}

In~\cite{invariant}, the $\I$-invariant is then defined as
\begin{equation*}
	\I(\A,\xi,\gamma)=\xi \circ \iota_* (\bar{\gamma}),
\end{equation*}
where $\iota_*$ is the map induced by the inclusion $\iota$ of the boundary $B_\A$ of a tubular neighbourhood $\Tub(\A)$ of $\A$ in $\CC\PP^2\setminus\A$, and where $\bar{\gamma}$ is a suitably chosen lift of the cycle $\gamma$ in $B_\A$. More precisely, this lift is a `nearby cycle' in the terminology of~\cite{invariant}, which roughly means that this cycle is contained in $B_{\C}\setminus \Tub(\D)\subset \CC\PP^2 \setminus \D$, where $\A=\C\cup \D$ with $\Irr(\C)=\{L\in\Irr(\A) \vert L\cap\gamma\neq\emptyset\}$ and $\Irr(\D)=\{L\in\Irr(\A) \vert L\cap\gamma = \emptyset\}$ -- see~\cite[Def.~2.11]{invariant} for a precise definition.

\begin{thm}
Let $\A=\C\cup\D$ be a line arrangement, and let $\gamma$ be a maximal cycle in $\Gamma_\C$.  
Let $\xi$ be a character on $\HH_1(\CC\PP^2\setminus\A)$ such that $(\A,\xi,\gamma)$ is an inner-cyclic arrangement. 
Then there is a nontrivial character $\xi_*$ on $\HH_1(\CC\PP^2\setminus\D) / \J_{(\C,\gamma)}$ induced by $\xi$ such that
	\begin{equation*}
		\I(\A,\xi,\gamma)=\xi_*(  \lk_\gamma(\C,\D) ).
	\end{equation*}
\end{thm}
\begin{proof}
Since $(\A,\xi,\gamma)$ is an inner-cyclic arrangement, we have that $\xi(m_L)=1$ for each line $L$ of $\A$ with corresponding component-vertex in $\gamma$  -- these correspond to the lines of $\C$ since $\gamma$ is maximal.
This shows that $\xi$ factors through the projection map $\HH_1(\CC\PP^2\setminus\A)\rightarrow \HH_1(\CC\PP^2\setminus \D)$. 
Furthermore, conditions (1) and (2-ii) of Definition~\ref{def:inner-cyclic} ensures that $\xi$ further factors to 
$\HH_1(\CC\PP^2\setminus\D) / \J_{(\C,\gamma)}$, thus providing the desired nontrivial character $\xi_*$.\footnote{This shows, in particular, that the quotient $\HH_1(\CC\PP^2\setminus\D) / \J_{(\C,\gamma)}$ is nontrivial. }
So we have 
$$ \I(\A,\xi,\gamma)=\xi \circ \iota_* (\bar{\gamma}) = \xi_* ([\bar{\gamma}]), $$
where $\bar{\gamma}$ is any lift of $\gamma$ which is a nearby cycle, and $[\bar{\gamma}]$ denotes its image in $\HH_1(\CC\PP^2\setminus\D) / \J_{(\C,\gamma)}$. 
On the other hand, we have by definition that 
$$  \lk_\gamma(\C,\D) = [\widetilde{\gamma}]\in  \HH_1(\CC\PP^2\setminus\D) / \J_{(\C,\gamma)}, $$
where $\widetilde{\gamma}$ is  a $\D$-avoiding lift of $\gamma$. 
The result then follows from the fact that the homology classes of the cycles $\bar{\gamma}$ and $\widetilde{\gamma}$ in $\HH_1(\CC\PP^2\setminus\D)$ can only differ by elements of $\J_{(\C,\gamma)}$, as follows from the definition of a nearby cycle given in \cite{invariant}.
\end{proof}

\section{Computations}\label{sec:computation}
	In this section, we describe two concrete methods for computating our invariant. The first one is topological and is based on a modification of the braid monodromy, and uses the usual linking number of links in the $3$-sphere. The second one is an algebraic method using the connected numbers introduced by Shirane in~\cite{shirane_connected}, and its relations with the linking invariant, observed in~\cite{GueShi}.

\subsection{Topological method}

We consider an algebraic curve $\C\cup \D$, such that $\C$ is rational, together with a cycle $\gamma$ in the intersection graph of $\C$.
\begin{de}
 A path $\widetilde{\gamma}$ in $\C$ is \emph{$\gamma$-admissible} if it is a lift of $\gamma$ 
 and if there is a generic projection $\pi: \CC\PP^2\setminus \{\ast\}\rightarrow \CC\PP^1$ such that $\pi(\widetilde{\gamma})$ has no self-intersection,  and $\left(\pi\inv\circ\pi(\widetilde{\gamma})\right) \cap \Sing(\D) =\emptyset$.
\end{de}

\noindent
Note that the latter condition can always be fulfilled, up to a small modification of $\pi$.

Let $\widetilde{\gamma}$ be a $\gamma$-admissible path in $\C$. For any point $p$ of $\pi(\widetilde{\gamma})$, we consider the fiber $F_p$ over~$p$. By the definition of $\gamma$-admissibility, the number of intersection points of $\D$ with $F_p$ equals the degree of $\D$ above all points of $\pi(\widetilde{\gamma})$. We denote by $L_\D$ the oriented link
\begin{equation}\label{eq:LD}
	L_\D=\left(\pi\inv\circ\pi(\widetilde{\gamma})\right)\cap \D \subset \CC\PP^2.
\end{equation}
Noting that $\widetilde{\gamma}$ and $L_\D$ do not intersect, we define $\widetilde{\gamma}\cup L_\D$. This link naturally sits in a copy of $S^3$, as follows. Let $D$ be the disc bounded by $\pi(\widetilde{\gamma})$ in $\CC\PP^1$. Pick a polydisc $\P$ of $\CC\PP^2$ such that $\pi(\P)= D$ and $\pi\inv(D)\cap (\C\cup\D) \subset \P$. 
By construction, the link $\widetilde{\gamma}\cup L_\D$ lies in the boundary of $\P$, which is homeomorphic to $S^3$. This construction yields the following. 
\begin{thm}\label{thm:topological_computation}
	Under the above assumptions, the homology class of the path $\widetilde{\gamma}$ in $\HH_1(\CC\PP^2\setminus\D)$
	 is given by:
	\begin{equation*}
		\rho\left(\sum_{c \in L_\D} \ell(c,\widetilde{\gamma}).m_c\right),
	\end{equation*}
	where $\ell$ is the usual linking number in $S^3$ and $\rho$ is the map induced by 
	the inclusion of $S^3\setminus L_\D$ in $\CC\PP^2\setminus \D$, which maps the meridian $m_c$ of each component $c$ of $L_\D$ to the meridian of the irreducible component of $\D$ containing $c$.
\end{thm}
\begin{rmk}
	If the curve $\C$ is not rational, we can still use the present method to compute the value of the generators of $\I_\C$ in $\HH_1(\CC\PP^2\setminus\D)$, and thus to compute the coset $[\gamma]\I_\C$.
\end{rmk}

This provides a computational formula for the linking invariant of $\C$ with $\D$ along $\gamma$ in terms of the usual linking number. See Section \ref{sec:quartic} for an application on a concrete example. 

\subsection{Algebraic method}

The second method of computation comes from the connected numbers introduced by Shirane in~\cite{shirane_connected}. This method applies when $\C$ is a nodal curve with $\Sing(\C)\cap\D=\emptyset$ (this implies that $\C\setminus\D$ is connected), and if $\HH_1(\CC\PP^2\setminus\D)/\J_\C\simeq\ZZ/m\ZZ$.

Let $\psi:X\rightarrow\CC\PP^2$ be a cyclic cover of degree $m$ branched over $\D$. The \emph{connected number} of $\C$ for $\psi$ is the number of connected components of $\psi^{-1}(\C\setminus\D)$) in $X$. Based on a previous version of the present paper, it has been proved in~\cite{GueShi} that if $\C$ and $\D$ are smooth curves then the connected number and the global linking of $\C$ and $\D$ are essentially equivalent. For the purpose of this paper, however, we will rather give the following statement. 
\begin{thm}\label{thm:algebraic_computation}
	Suppose that $\C$ is a nodal curve such that $\Sing(\C)\cap\D=\emptyset$, and that $\HH_1(\CC\PP^2\setminus\D)/\J_\C\simeq\ZZ/m\ZZ$.
	Then the global linking of $\C$ and $\D$ is the unique subgroup of $\HH_1(\CC\PP^2\setminus\D)/\J_\C$ of index $\deg(\D)/\mu$, where $\mu$ is the minimal degree of a plane curve $E$ 
	such that 
$\C\cap\D=\C\cap E$ and for each point $P\in\C\cap\D$, we have $(\deg(\D)/\mu).I_P(\C,E) = I_P(\C,\D)$.
\end{thm}

\begin{proof}[Sketch of proof]
	By the proof of~\cite[Theorem~2.5]{GueShi}, we know that the index of $\L(\C,\D)$ in $\HH_1(\CC\PP^2\setminus\D)/\J_\C$ is equal to the connected number of $\C$ for $\psi$. Since $\HH_1(\CC\PP^2\setminus\D)/\J_\C$ is cyclic then each of its subgroup is determined by its index. The result is thus a consequence of~\cite[Corollary~2.5]{shirane_connected}.
\end{proof}

\section{Applications}\label{sec:application}
	In this section, we use our linking invariants (both the oriented and global versions) to distinguish two types of Zariski pairs. 

\subsection{$3$-Artal curves}

As an application of the algebraic method of computation of the linking invariant, we propose to distinguish the Zariski pair found by Artal in~\cite{artal}. These curves are formed by a smooth cubic $\C$ and three inflexional tangent lines. The geometry of the 9 inflexion points of a cubic is well known; the collinearity relations are the same as in $\FF_3^2$, the plane over the finite field of 3 elements. We consider $P_1,\dots,P_4$ four inflexion points of $\C$ such that $P_1,P_2,P_3$ are collinear and $P_1,P_2,P_4$ are not. Set $\A_1=L_1 \cup L_2 \cup L_{3}$ and $\A_2=L_1 \cup L_2 \cup L_{4}$, where $L_i$ denotes the inflexional tangent line at $P_i$.

\begin{thm}
	The global linking of the line arrangement $\A_i$ with the cubic $\C$ is:
	\begin{equation*}
		\L(\A_i,\C)=
		\left\{
			\begin{array}{ccl}
				\{0\} & \quad & \text{if } i=1,\\
				\{0,1,2\} & \quad & \text{if } i=2.
			\end{array}
		\right.
	\end{equation*}
\end{thm}

\begin{proof}
	Since the cubic $\C$ is smooth, we have $\HH_1(\CC\PP^2\setminus\C)\simeq\ZZ/3\ZZ$. Furthermore, Proposition~\ref{propo:indeterminacy} implies that $\J_{\A_i}=\langle 3m \rangle$, where $m$ is a meridian of the cubic. So the quotient $\HH_1(\CC\PP^2\setminus\C)/\J_{\A_i}$ is also isomorphic to $\ZZ/3\ZZ$. 
	
	We first compute the global linking $\L(\A_1,\C)$. By construction, there is a line $E$ (\emph{i.e.} an algebraic curve of degree~1) passing through $P_1,P_2$ and $P_3$. By Bezout theorem, we have $I_{P_j}(\A_1,E)=1$ for $j\in\{1,2,3\}$, and the following equality holds for $j\in\{1,2,3\}$ 
	\begin{equation*}
		\frac{\deg(\C)}{\deg(E)}.I_{P_j}(\A_1,E) = I_{P_j}(\A_1,\C).
	\end{equation*}
	Thus by Theorem~\ref{thm:algebraic_computation} the index of $\L(\A_1,\C)$ in $\HH_1(\CC\PP^2\setminus\C)/\J_{\A_1}$ is 3.
	
	Let us now turn to $\L(\A_2,\C)$, and look for the minimal degree curve $E$  passing through the points $P_1,P_2$ and $P_4$ and satisfying $\frac{\deg(\C)}{\deg(E)}.I_{P_j}(\A_2,E) = I_{P_j}(\A_2,\C)$ for $j\in\{1,2,4\}$. By construction, no line $E$ satisfies these conditions. Similar considerations as above show that no conic can verify these conditions either\footnote{This also follows from Theorem~\ref{thm:algebraic_computation}, since the existence of such a conic would imply that $\L(\A_1,\C)$ is an index~2 subgroup of $\ZZ/3\ZZ$.}. But taking $E$ to be the cubic $\C$ obviously works, and it follows by Theorem~\ref{thm:algebraic_computation} that the index of $\L(\A_2,\C)$ in $\HH_1(\CC\PP^2\setminus\C)/\J_{\A_2}$ is 1.
\end{proof}

\begin{cor}
	The curves $\C\cup\A_1$ and $\C\cup\A_2$ form a Zariski pair.
\end{cor}

\subsection{The quartic and its bitangents}\label{sec:quartic}

As an application of the topological method, we will distinguish a Zariski pair formed by a quartic and 3 bitangents. Let $\Q$ be the Klein quartic defined by $x^3y+y^3z+z^3x=0$. The full list of its 28 bitangents is given in~\cite{Shioda}. We will consider here only four of them. Let $\zeta$ be a primitive $7$th root of unity, and define the real numbers $\epsilon_i=\zeta^i+\zeta^{-i}$, for $i\in\{1,2,3\}$. We consider the following bitangents:
\begin{equation*}
	\begin{array}{lll}
		L_1:\ x+y+z=0, & \quad & L_2:\ x+\epsilon_3^{-2}y+\epsilon_1^{2}z=0, \\
		L_3:\ x+\zeta^3\epsilon_3^{-2}y+\zeta\epsilon1^{2}z=0, & \quad & L_4:\ x+\zeta^2\epsilon_2^{-2}y+\zeta^3\epsilon_3^{2}z=0. 
	\end{array}
\end{equation*}
Let $\A_i=L_1 \cup L_2 \cup L_{2+i}$, for $i\in\{1,2\}$. We will compute the linking of $\A_i$ with $\Q$ along $\gamma_i$, where $\gamma_i$ is a cycle generating $\HH_1(\A_i)\simeq \mathbb{Z}$. Since $\HH_1(L_j)=0$ for all $j$, we have $\I_{\A_i}=0$ for $i=1,2$. Furthermore, $\Q$ is smooth so $\HH_1(\CC\PP^2\setminus\Q)\simeq\ZZ/4\ZZ$. By Proposition~\ref{propo:indeterminacy}, we have that $\J_{\A_i}=\langle 2m \rangle$, where $m$ is a meridian of $\Q$. So we have 
\begin{equation*}
	\HH_1(\CC\PP^2\setminus\Q)/\J_{\A_i}\simeq\ZZ/2\ZZ.
\end{equation*}

By Theorem \ref{thm:topological_computation}, the linking invariant $\lk_{\gamma_i}(\A_i,\Q)$ is given by the linking numbers of the cycle $\gamma_i$ with each component of the link $L_\Q$ defined in (\ref{eq:LD}). The cycle $\gamma_i$ can be viewed as the triangle formed by the three lines, and in order to determine the desired linking numbers, we can glue a $2$-cell $c_i$ along $\gamma_i$ and compute its intersection with the link $L_\Q$.

\begin{figure}[ht!]
	\centering
		\includegraphics[width=6cm]{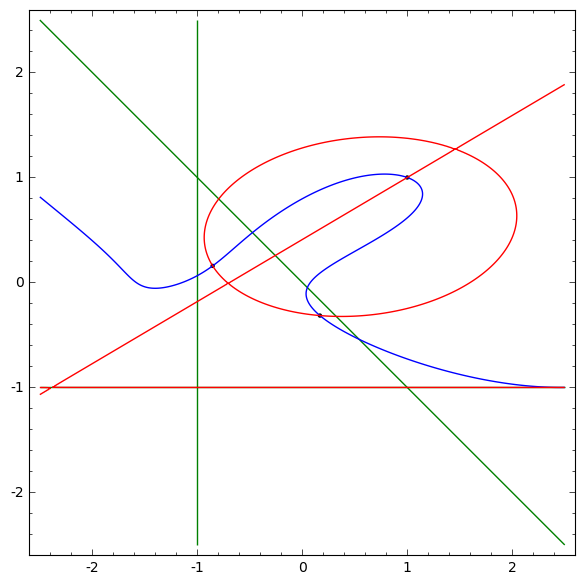}
		\hspace{1cm}
		\includegraphics[width=6cm]{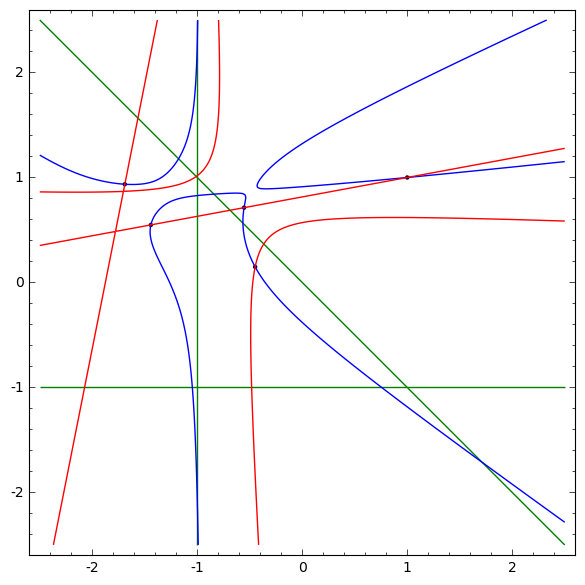}\\
		$i=1$ \hspace{6cm} $i=2$
	\caption{Real part of $\Q\cup\A_i$ \label{fig:Q} ($\A_i$ is in green, $\Q_{Re}$  in blue and $\Q_{Im}$ in red) }
\end{figure}

In order to simplify the computation, we apply a linear change of variable $P_i$ such that the line arrangement $\A_i$ is given by the real equation $(x+y)(y+z)(z+x)=0$. With this modification, we can take as the $2$-cell $c_i$ the interior of the real triangle formed by $\A_1$ (resp. $\A_2$) in the chart $z=1$ (resp. $y=1$). After the change of variable $P_i$, the quartic $\Q$ can be written as the sum of its real and complex part $\Q=\Q_{Re}+i\Q_{Im}$, where $\Q_{Re}$ and $\Q_{Im}$ are quartic with real coefficients. The real points $[x_0:y_0:z_0]$ of $\Q$ are those which verify $\Q_{Re}(x_0,y_0,z_0)=0$ and $\Q_{Im}(x_0,y_0,z_0)=0$. Plotting the arrangement $\A_i$ and the two curves $\Q_{Re}$ and $\Q_{Im}$ (see Figure~\ref{fig:Q}), we compute the linking of $\gamma_i$ with $L_\Q$ counting the parity of the intersection points $\Q_{Re}$ and $\Q_{Im}$ in the triangle defined by $\A_i$ and we get that the value of $\gamma_i$ in $\HH_1(\CC\PP^2\setminus\Q)/\J_{\A_i}$ is
\begin{equation*}
	\rho(\gamma_i)=
	\left\{
		\begin{array}{lll}
			0 &\quad& \text{if } i=1,\\
			1 &\quad& \text{if } i=2.
		\end{array}
	\right.
\end{equation*} 

From these computations and using Remark~\ref{rem:non-oriented}, we have the following theorem and its corollary.
\begin{thm}
	The linking of the line arrangement $\A_i$ with the quartic $\Q$ along $\gamma_i$ is 
	\begin{equation*}
		\lk_{\gamma_i}(\A_i,\Q)=
		\left\{
			\begin{array}{cll}
				\{0\} &\quad& \text{if } i=1,\\
				\{1\} &\quad& \text{if } i=2.
			\end{array}
		\right.
	\end{equation*} 
\end{thm}
\begin{cor}
	The curves $\Q\cup\A_1$ and $\Q\cup\A_2$ form a Zariski pair.
\end{cor}
\begin{rmk}
	Since $\gamma_i$ generates $\HH_1(\A_i)$, we can also compute the global linking of $\A_i$ with $\Q$. It is given by:
	\begin{equation*}
		\L(\A_i,\Q)=
		\left\{
			\begin{array}{cll}
				\{0\} &\quad& \text{if } i=1,\\
				\{0,1\} &\quad& \text{if } i=2.
			\end{array}
		\right.
	\end{equation*} 
\end{rmk}

\begin{rmk}
	This result is in adequation with the computation of the connected numbers of these curves made in~\cite{BTY}.
\end{rmk}

\subsection*{Acknowledgement}
The authors would like to thanks E.~Artal for seminal discussions at an early stage of this work.
They also thank S.~Bannai, V.~Florens, P.~Popescu-Pampu, T.~Shirane and H.~Tokunaga for their questions and remarks. \\


\bibliographystyle{abbrv}
\bibliography{paper}

\end{document}